\documentclass[11pt]{article}

\usepackage[T1]{fontenc}
\usepackage{amssymb}
\usepackage{amsmath}
\usepackage[title, toc]{appendix}
\usepackage{amsthm}
\usepackage{tikz}
\usepackage{tikz-cd}  
\usepackage{enumerate} 
\usepackage{enumitem}  
\usepackage{pdfpages}  
\usepackage{hyperref}
\usepackage{wasysym}  
\usepackage{mathtools} 
\usepackage{stmaryrd} 
\usepackage{graphicx}
\usepackage{xcolor}
\usepackage{mathrsfs}
\usepackage{titlesec}
\usepackage{todonotes} 
\usepackage{xcolor}
\hypersetup{
    colorlinks,
    linkcolor={blue!50!black},
    citecolor={blue!50!black},
    urlcolor={blue!80!black}
}
\DeclareSymbolFont{extraitalic} {U}{zavm}{m}{it}
\DeclareMathSymbol{\stigma}{\mathord}{extraitalic}{168}

\newtheorem*{proposition*}{Proposition}
\newtheorem*{theorem*}{Theorem}
\newtheorem{theorem}{Theorem}

\newtheorem{lemma}[theorem]{Lemma}
\theoremstyle{remark}
\newtheoremstyle{boldremark}
    {\dimexpr\topsep/2\relax} 
    {\dimexpr\topsep/2\relax} 
    {}          
    {}          
    {\bfseries} 
    {.}         
    {.5em}      
    {}          
\theoremstyle{boldremark}
\newtheorem*{conjecture*}{Conjecture}
\makeatletter
\newcommand*\leftdash{\rotatebox[origin=c]{-45}{$\dabar@\dabar@\dabar@$}}
\newcommand*\rightdash{\rotatebox[origin=c]{45}{$\dabar@\dabar@\dabar@$}}
\makeatother

\newcommand{\Qlbar}{\overline{\mathbf{Q}}_{\ell}}

\newcommand{\kk}{\mathbf{k}}

\newcommand{\GL}{\textrm{GL}}

\newcommand{\h}{\mathfrak{h}}

\newcommand{\g}{\mathfrak{g}}

\newcommand{\levi}{\mathfrak{l}}

\newcommand{\gl}{\mathfrak{gl}}

\newcommand{\F}{\mathcal{F}}
\newcommand{\G}{\mathcal{G}}

\newcommand{\cc}{\underline{\Qlbar}}
\newcommand{\cF}{\Qlbar}

\newcommand{\gspr}{\mathfrak{gspr}}
\newcommand{\spr}{\mathfrak{spr}}

\newcommand{\FT}{\operatorname{FT}}

\newcommand{\Ad}{\operatorname{Ad}}

\titleformat{\section}[runin]
{\normalfont\bfseries}{\thesection}{1em}{}
\titleformat{\subsection}[runin]
{\normalfont\bfseries}{\thesubsection}{1em}{}
\titleformat{\subsubsection}[runin]
{\normalfont\bfseries}{\thesubsection}{1em}{}
 \newcommand{\Addresses}{{
  \bigskip
  \footnotesize

  R.~Bezrukavnikov, \textsc{Department of Mathematics, Massachusetts Institute of Technology,
    Cambridge, Massachusetts}\par\nopagebreak
  \textit{E-mail address}: \texttt{bezrukav@math.mit.edu}

  \medskip

  K.~Tolmachov, \textsc{School of Mathematics, University of Edinburgh,
    Edinburgh, United Kingdom}\par\nopagebreak
  \textit{E-mail address}: \texttt{tolmak@khtos.com}
 }}
 \title{Exterior powers of a parabolic Springer sheaf on a Lie algebra}

\author{Roman Bezrukavnikov, Kostiantyn Tolmachov} \date{}
\begin{document}
\maketitle \abstract{
 We compute the exterior powers, with respect to the additive convolution on the general linear Lie algebra, of a parabolic Springer sheaf corresponding to a maximal parabolic subgroup of type (1, n -- 1). They turn out to be isomorphic to the semisimple perverse sheaves attached by the Springer correspondence to the exterior powers of the permutation representation of the symmetric group. 
}
\section{Basic notations.}
Fix primes $\ell\neq p$. For a stack $X$ defined over an algebraically closed field $\kk$ of characteristic $p$, let $D^b(X)$ be the bounded derived category of $\Qlbar$-sheaves with constructible co\-ho\-mo\-lo\-gy. All stacks in question will be quotient stacks for an action of an algebraic group on a nice scheme, and we use \cite{bernsteinEquivariantSheavesFunctors2006}, \cite{laszloSixOperationsSheaves2008a} for the theory of $D^b(X)$.  Let $P(X)$ be the abelian category of perverse sheaves on $X$ with respect to the middle perversity. Let $\cc_{X}$ stand for the constant sheaf on $X$.

All our results can be also stated, with standard modifications, for algebraic varieties and stacks over $\mathbb{C}$, with $D^b(X)$ being the bounded derived category of $\mathbb{Q}$-sheaves in analytic topology with algebraically-constructible cohomology. Formulation of the theory of Fourier-Deligne transform in Section \ref{sec:four-deligne-transf} should be then adjusted in a standard way.

$G$ will stand for a reductive group defined over $\kk$, and $\g$ for its Lie algebra. If $H, \h$ are an algebraic group and its Lie algebra, respectively, we will write $\h/H$ for the adjoint quotient stack.

\section{Introduction.}
The category $D^b(\g)$ has a symmetric monoidal structure coming from the additive convolution. Namely, let \[a:\g\times \g \to \g, \pi_{1,2}:\g\times\g \to \g\] be the addition map and two projections, respectively. For $\F, \G \in D^b(\g)$, write
\begin{equation}
  \label{eq:1}
  \F \star \G = m_!(\pi_1^*\F \otimes \pi_2^*\G).
\end{equation}

The map $a$ is $G$-equivariant with respect to the diagonal adjoint action of $G$ on $\g \times \g$ and adjoint action on $\g$. Thus \eqref{eq:1} equips the category $D^b(\g/G)$ with the symmetric monoidal structure. For an object $\F \in D^b(\g/G),$ we write $\wedge^k_{+}\F$ for its $k$th exterior power with respect to this structure.

Let now $\g = \gl_n$ be the general linear Lie algebra over $\kk$. 

Let $\spr(M)$ be the semisimple sheaf attached by the Springer correspondence to a representation $M$ of the symmetric group $S_n$ -- see the next section for the reminder of definitions. Let $V$ stand for the permutation representation of $S_n$.

Main goal of this note is to give a simple proof of the following result.
\begin{theorem}
\label{sec:stat-main-results-1}
  Exterior powers of $\spr(V)$ with respect to the symmetric monoidal structure $\star$ on $D^b(\gl_n/G)$ satisfy
  \[\wedge^k_+\spr(V) \simeq \spr(\wedge^k V).\]
  In particular, $\wedge^{k}_+\spr(V) = 0$ for $k > n$. 
\end{theorem}

This result, with a more complicated proof, appeared in the se\-cond author's thesis \cite{tolmachovFunctorAffineFinite2018}. This is a linearization of a similar, but more technically involved, statement for a parabolic Springer sheaf on the group $\GL_n$ and exterior powers with respect to the multiplicative convolution, see loc.cit. and \cite{tolmachovkostiantynLinearStructureFinite}. One surprising outcome of these computations is that exterior powers of the parabolic Springer sheaf with respect to the convolution operation remain perverse and semisimple.

We believe that linearized result, together with a simple proof presented, is of independent interest. Note that in \cite{hhh} it is proved that sheaves attached by the Springer correspondence to the exterior powers of the representation of the Weyl group on a Cartan subalgebra appear as perverse cohomology of the $G$-averaging of a Whittaker sheaf on a maximal unipotent subgroup, for any reductive group $G$. We don't know of any connection between the methods of loc. cit. and the present note, or if the result of the present note can be extended to a more general $\g$.

\section{Springer theory.} Let for now $G$ be any reductive group. Let $\g^{rss}$ stand for the subset of regular semisimple elements in $\g$. We recall some basic facts from Springer theory. See \cite{AST_1983__101-102__23_0} and references therein. 

For a parabolic $P \subset G$, let $U_P$ be its unipotent radical. Let $\mathfrak{p}, \mathfrak{n}_P$ be the Lie algebras of $P,U_P$. Let \[\tilde{\mathcal{N}}_{P} = G \times^P \mathfrak{n}_P, \tilde{\g}_P = G \times^P\mathfrak{p}\] be the parabolic Springer and Grothendieck-Springer varieties, and $q: \tilde{\mathcal{N}}_{P} \to {\g}, q': \tilde{\g}_P \to \g$ stand for maps given by $(g,x) \mapsto \Ad_g x$. Denote by \[\mathfrak{spr}_P:= q_*\cc_{\tilde{\mathcal{N}}_{\g}}[\dim \mathcal{N}_{P}], \mathfrak{gspr}_P:=q'_*\cc_{\tilde{\g}_P}[\dim \g],\] the parabolic Springer and Grothendieck-Springer sheaves. They are perverse and semisimple. $\mathfrak{spr}_P$ is supported on the nilpotent variety $\mathcal{N}_\g$ of $\g$.  $\mathfrak{gspr}_P$ is a local system over the open subset $\g^{rss}$, and is the intermediate extension of its restriction to $\g^{rss}$.

All maps involved are $G$-equivariant with respect to standard actions, and we will denote the perverse parabolic Springer and Grothendieck-Springer sheaves on $\g/G$ in the same way.

Let $B$ be a Borel subgroup, and let $W$ be the Weyl group of $G$. We have
\[
  \operatorname{End}(\spr_B) \simeq \operatorname{End}(\gspr_B) \simeq \cF[W].
\]
There are two standard ways to choose an action of $W$ on $\spr_B, \gspr_B$, that differ by a sign representation of $W$. We fix our action by fixing the $W$-invariant summand:  \[\spr_B^W = \iota_*\cc_{\operatorname{pt}}, \gspr_B^W = \cc_\g[\dim \g],\] where $\iota: \operatorname{pt} \to \g$ is the embedding of $0$. 

For a representation $V$ of $W$ we write \[\spr(V) = V\otimes_{\cF[W]}\spr_B, \gspr(V) = V\otimes_{\cF[W]}\gspr_B.\]

From now on, let $\g = \gl_n$ be the Lie algebra of the reductive group $G = \GL_n$. Let $V$ be an $n$-dimensional vector space. We identify $G$ with $\GL(V)$.

Fix a flag $V_1 \subset V_2 \subset \dots \subset V_{n-1} \subset V_n = V$ with
$\dim V_k = k$, and let let $P_k$ be the parabolic subgroup of $\GL_n$ preserving
$V_k$. Let $L_k$ be its Levi subroup, and $\levi_k$ the Lie algebra of $L_k$. 

Fix a basis of $V$ compatible with the flag above. This equips $V$ with the permutation representation of $S_n \simeq W$. We have
\[
  \spr_{P_1} = \spr(V), \gspr_{P_1} = \gspr(V).
\]

\section{Fourier-Deligne transform.}\label{sec:four-deligne-transf}
The categories $D^b(\g), D^b(\g/G)$ are also equipped with the standard monoidal structure of tensor product $\otimes_{\cF}$.  

Recall the Fourier-Deligne transform functor
    \[
      \operatorname{FT}:D^b(\g/G) \to D^b(\g/G),
    \]
    from  \cite{katzTransformationFourierMajoration1985}, \cite{brylinskiTransformationsCanoniquesDualite1986}. $\operatorname{FT}$ is an equivalence, intertwines the convolution monoidal structure with the shifted tensor product $(-\otimes_{\cF}-)[-\dim\g]$, and we have
    \begin{equation}
      \label{eq:2}
      \operatorname{FT}(\spr(V)) \simeq \gspr(V).
    \end{equation}
  For an object $\F \in D^b(\g/G),$ we write $\wedge^k_{\otimes}\F$ for its $k$th exterior power with respect to the shifted tensor product $(-\otimes_{\cF}-)[-\dim\g]$.  Theorem \ref{sec:stat-main-results-1} is therefore equivalent to 
\begin{theorem}
\label{sec:stat-main-results-2}
  Exterior powers of $\gspr(V)$ with respect to the symmetric monoidal structure $\otimes_{\cF}$ on $D^b(\gl_n/G)$ satisfy
  \[\wedge^k_{\otimes}\gspr(V) \simeq \gspr(\wedge^k V).\]
  In particular, $\wedge^{k}_{\otimes}\gspr(V) = 0$ for $k > n$. 
\end{theorem}
Let us first show that $\wedge^{k}_{\otimes}\gspr(V) = 0$ for $k > n$. Recall the map $q':\gspr_{P_1} \to \g$.  
It is easy to see that the fiber $\gspr_{P_1,x}$ of $q'$ over any $x \in \g$ is a union of projective spaces having the total dimension of cohomology bounded by $n$. Let $i_x:\operatorname{pt} \to \g$ be the inclusion of a point $x \in \g$. We have
\[
  i_x^*\wedge^{k}_{\otimes}\gspr(V) = \wedge^{k}i_x^*\gspr(V) = \wedge^k \operatorname{H}_c^{\bullet}(\gspr_{P_1,x}) = 0,
\]
and the claim follows.

Note now that the statement of Theorem \ref{sec:stat-main-results-2} is
readily checked to be true over the open subset $\g^{rss}$. It is thus enough
to show that $\wedge^k_{\otimes}\gspr(V)$ is perverse and is given by the
intermediate extension of its restriction to $\g^{rss}$. We will now assume
that the Theorem is proved for $G = \GL_m, m < n$. 

\section{\'{E}tale neighbourhoods in $\g$.} 
Let $P$ be an arbitrary parabolic subgroup of $G$, let $L$ be its Levi subgroup and $\levi$ be the Lie algebra of $L$. Recall that the element $x \in \g$ with semisimple part $x_s$ is called 
$L$-regular if the centralizer $C_G(x_s)$ is conjugate to a subgroup of $L$. Write $\levi^{r}$ for the set of $L$-regular elements of $\g$ lying inside $\levi$. We will use the following
\begin{lemma}[Proposition 2.11 in \cite{gunninghamGeneralizedSpringerTheory2018}]
  \label{lemma-etale} The natural map $\mathfrak{l}^r/L \to \g/G$ is \'{e}tale, with image consisting of the set of $L$-regular elements in $\g$. 
\end{lemma}

For $m = 1, \dots, n-1$ let $j_m: \mathfrak{l}_m^r/L_m \to \g/G$ be the natural map. It follows from Lemma \ref{lemma-etale} that $j_r$ together give an 
\'{e}tale open cover of $\g\backslash \mathcal{N}_{\g}$. Note that $\mathfrak{l}_m\simeq \mathfrak{gl}_m \times \mathfrak{gl}_{n-m}$, and it is easy to check
that $j_m^*\gspr(V) \simeq \gspr(\operatorname{Res}^{S_n}_{S_m \times S_{n-m}}V)|_{\mathfrak{l}_m^r}$. By our inductive assumption, it follows that $j_m^*\wedge^k_{\otimes}\gspr(V)$
is given by the intermediate extension of its restriction to $\g^{rss}$.

To finish the proof, it is now enough to show that perverse cohomology of $\wedge^k_\otimes\gspr(V)$ have no perverse constitutents supported on ${\mathcal{N}_\g}$.

\section{Rank estimate.} Any perverse sheaf on $\mathcal{N}_{\gl_n}$ is a
direct sum of sheaves of the form $\spr(M), M \in \operatorname{Rep}(W)$. Since $\FT(\spr(M))
\simeq \gspr(M)$ and $\gspr(M)$ has full support, it is enough to show that no perverse constituents of
$\FT(\wedge^k_\otimes\gspr(V))\simeq\wedge^k_{+}\spr(V)$ have full support. 

Note that $\spr(V) \simeq \spr(\operatorname{triv})\oplus\spr(R)$, where $R$
stands for the reflection representation of $W$. Since
$\spr(\operatorname{triv}) \simeq \iota_*\cc_{\operatorname{pt}}$ is the monoidal unit of the additive convolution,
we have 
\[
  \wedge^k_+\spr(V)\simeq \wedge^k_+\spr(R)\oplus\wedge^{k-1}_+\spr(R).
\]
Since $\wedge^{n+1}_+\spr(V) = 0$, we have $\wedge^n_+\spr(R) = 0$, and so it is enough to show that that $\wedge^k_+\spr(R)$ has no perverse constituents
with full support for $k = 1, \dots, n-1$. 

Note that, since all elements of $\mathfrak{u}_{P_1}$ have rank bounded by $1$, the sum of $n-1$ such elements has rank bounded by $n-1$. It follows that $\wedge^{k}_+\spr(R)$ is supported on degenerate matrices for $k < n$, and hence does not have full support. This finishes the proof of Theorems \ref{sec:stat-main-results-1} and \ref{sec:stat-main-results-2}. 
\bibliography{bibliography} \bibliographystyle{alpha}
 \Addresses
\end{document}